# A Simplification in the proof presented for non existence of periodic solutions in time invariant fractional order systems

Mohammad Saleh Tavazoei[a] and Mohammad Haeri[b]

[a] Electrical Engineering Department, Sharif University of Technology (tavazoei@sharif.edu)
[b] Advanced Control System Lab, Electrical Engineering Department, Sharif University of Technology (haeri@sina.sharif.edu)

**Abstract**

In this note, a short-cut is proposed to shorten the proof which has been previously presented for non existence of periodic solutions in time invariant fractional order systems.

In Tavazoei & Haeri (2009), a proof for non existence of periodic solutions in time invariant fractional order systems has been presented. The aim of this note is proposing a short-cut to simplify the mentioned proof. To this end, let the following Lemma be proposed after proving Lemmas 1 and 2 in Tavazoei & Haeri (2009) (It should be noted that all notations and assumptions in this note are similar as those of Tavazoei & Haeri (2009)).

**Lemma**: If equation (10) in Tavazoei & Haeri (2009) is true for all $p \in \mathbb{N}$, then

$$\int_0^T \tau^m \tilde{g}_i(\tau) d\tau = 0, \quad i=1,2,\cdots,n, \tag{1}$$

where $m \in \mathbb{N}$.

**Proof:** If $h_i(t)$ is a positive function in the range $[0,T]$ and $c$ is a constant greater than $T^m$, then it can be easily shown that

$$m_i \bar{S}_i^+ - M_i \bar{S}_i^-$$
$$\leq \int_0^T h_i(\tau)(c-\tau^m)\tilde{g}_i(\tau)d\tau \leq \qquad , \tag{2}$$
$$M_i \bar{S}_i^+ - m_i \bar{S}_i^-, \quad i=1,2,\cdots,n$$

where

$$m_i = \min\{h_i(t), \ 0 \leq t \leq T\}, \quad i=1,2,\cdots,n, \tag{3}$$

$$M_i = \max\{h_i(t), \ 0 \leq t \leq T\}, \quad i=1,2,\cdots,n, \tag{4}$$

and $\bar{S}_i^+$ and $\bar{S}_i^-$ indicate, respectively, the area between the curve $(c-\tau^m)\tilde{g}_i(\tau)$ and $t$-axis located above the axis and the area between the curve $(c-\tau^m)\tilde{g}_i(\tau)$ and $t$-axis located below the axis. If $h_i(t) = (pT-t)^{\alpha_i-1}/(c-t^m)$, we have $m_i = (pT)^{\alpha_i-1}/c$ and $M_i = (pT-T)^{\alpha_i-1}/(c-T^m)$. In this case, (2) can be rewritten in the following form

$$\frac{(pT)^{\alpha_i-1}}{c}S_i^+ - \frac{(pT-T)^{\alpha_i-1}}{c-T^m}S_i^-$$
$$\leq \int_0^T (pT-\tau)^{\alpha_i-1}\tilde{g}_i(\tau)d\tau \leq \qquad . \tag{5}$$
$$\frac{(pT-T)^{\alpha_i-1}}{c-T^m}S_i^+ - \frac{(pT)^{\alpha_i-1}}{c}S_i^-$$

According to (5) and Lemma 1 in Tavazoei & Haeri (2009), it is deduced that

$$\frac{(pT)^{\alpha_i-1}}{c}S_i^+ - \frac{(pT-T)^{\alpha_i-1}}{c-T^m}S_i^- \leq 0$$
$$\frac{(pT-T)^{\alpha_i-1}}{c-T^m}S_i^+ - \frac{(pT)^{\alpha_i-1}}{c}S_i^- \geq 0 \tag{6}$$

which is equivalent to

$$\frac{(pT)^{\alpha_i-1}}{(pT-T)^{\alpha_i-1}}S_i^+ - \frac{c}{c-T^m}S_i^- \leq 0$$
$$\frac{(pT-T)^{\alpha_i-1}}{(pT)^{\alpha_i-1}}S_i^+ - \frac{c-T^m}{c}S_i^- \geq 0 \tag{7}$$

Since $\lim_{p \to \infty}(pT)^{\alpha_i-1}/(pT-T)^{\alpha_i-1}=1$, from (7) it is resulted that

$$S_i^+ - \frac{c}{c-T^m}S_i^- \leq 0 \leq S_i^+ - \frac{c-T^m}{c}S_i^-, \tag{8}$$

We know that $\int_0^T (c-\tau^m)\tilde{g}_i(\tau)d\tau = S_i^+ - S_i^-$. Hence, (8) results in

$$\left(\frac{c-T^m}{c}-1\right)S_i^- \leq \int_0^T (c-\tau^m)\tilde{g}_i(\tau)d\tau \leq \left(\frac{c}{c-T^m}-1\right)S_i^- \quad (9)$$

From (9), it is deduced that

$$\lim_{c\to\infty} \int_0^T (c-\tau^m)\tilde{g}_i(\tau)d\tau = 0. \quad (10)$$

From (10) and Lemma 2 in Tavazoei & Haeri (2009), (1) is obtained. ∎

A straightforward result of the proved lemma in this note and Lemma 2 in Tavazoei & Haeri (2009) is the following equalities

$$\int_0^T (\tau-t_0)^m \tilde{g}_i(\tau)d\tau = 0, \quad i=1,2,\cdots,n, \quad (11)$$

for every $t_0 \in (0,T)$ and $m \in \mathbb{N}\cup\{0\}$. Now, by using (11) and equalities (48)-(51) in Tavazoei & Haeri (2009) the main theorem in the mentioned paper is proved. Clearly, such an approach shortens the intended proof.

It is worthy mentioning that the idea used in this note can be also applied to more simplify the proof of Theorem 1 in Tavazoei (2010). In the mentioned theorem, it has been proved that fractional-order derivatives of a periodic function with a specific period can not be a periodic function with the same period.